\documentclass[twocolumn,10pt]{IEEEtran}
\usepackage{amsmath, amsthm, amssymb}
\usepackage{graphicx}
\usepackage{epstopdf}
\usepackage{overpic}
\usepackage{cancel}
\usepackage{rotating}
\usepackage{url}
\usepackage{caption}
\usepackage{color}
\usepackage{rotating}
\usepackage{multirow}
\usepackage{wrapfig}

\usepackage[bottom,flushmargin,hang,multiple]{footmisc}
\usepackage{lipsum}
\newcommand\blfootnote[1]{%
  \begingroup
  \renewcommand\thefootnote{}\footnote{#1}%
  \addtocounter{footnote}{-1}%
  \endgroup
}

\newcommand{\bA}{\mathbf{A}}
\newcommand{\bB}{\mathbf{B}}

\newcommand{\bx}{\mathbf{x}}

\newcommand{\bu}{\mathbf{u}}

\newcommand{\bU}{\mathbf{U}}
\newcommand{\bSigma}{\mathbf{\Sigma}}
\newcommand{\bUp}{\mathbf{\Upsilon}}
\newcommand{\bV}{\mathbf{V}}
\newcommand{\bX}{\mathbf{X}}

\definecolor{blue}{rgb}{0,0,1}
\definecolor{darkgreen}{rgb}{0,0.5,0}
\definecolor{red}{rgb}{1,0,0}
\definecolor{teal}{rgb}{0,0.5,0.7}

\graphicspath{{Figures/}}

\title{Dynamic mode decomposition with control}
\author{Joshua L. Proctor$^{1*}$, Steven L. Brunton$^2$, J. Nathan Kutz$^2$\\
\small{$^1$Institute of Disease Modeling Bellevue, WA 98004, United States}\\ 
\small{$^2$Applied Mathematics, University of Washington, Seattle, WA 98195, United States}\\

\normalsize{DRAFT: last updated \today}
}
\date{last compiled \today}                                           

\begin{document}
\maketitle
\blfootnote{$^*$ Corresponding author. Tel.: +1 509 868 5696.\\ {\indent\emph{E-mail address:} joproctor@intven.com (J.L. Proctor).}}
\begin{abstract}
We develop a new method which extends Dynamic Mode Decomposition (DMD) to incorporate the effect of control to extract low-order models from high-dimensional, complex systems.   
DMD finds spatial-temporal coherent modes, connects local-linear analysis to nonlinear operator theory, and provides an equation-free architecture which is compatible with compressive sensing.
In actuated systems, DMD is incapable of producing an input-output model; moreover, the dynamics and the modes will be corrupted by external forcing.
 Our new method, Dynamic Mode Decomposition with control (DMDc), capitalizes on all of the advantages of DMD and provides the additional innovation of being able to disambiguate between the underlying dynamics and the effects of actuation, resulting in accurate input-output models.
The method is data-driven in that it does not require knowledge of the underlying governing equations, only snapshots of state and actuation data from historical, experimental, or black-box simulations.   
We demonstrate the method on high-dimensional dynamical systems, including a model with relevance to the analysis of infectious disease data with mass vaccination (actuation).  

\end{abstract}

\section{Introduction}

We introduce the method of Dynamic Mode Decomposition with control (DMDc) to analyze observational data arising from complex, high-dimensional systems that exhibit dynamics and require control.  By utilizing both measurements of the system and the applied external control, the underlying, unforced dynamics can be extracted and specified in an {\em equation-free} manner, i.e. the underlying equations of motion do not have to be known.   In addition, a description of how the control inputs affect the system are also discovered and characterized.  With a quantitative understanding of the input-output characteristics, a reduced-order-model can be generated for both prediction and design of controllers for high-dimensional, complex systems.

Controlling high-dimensional systems remains an extremely challenging task as many control strategies do not scale well with the dimension of the system.  
In particular, controllers developed on a full system may be computationally prohibitive to implement, introducing unacceptably large latencies.  Moreover, many control laws are determined by solving a large Riccati equation ($\mathcal{H}_2$) or through an iterative procedure ($\mathcal{H}_{\infty}$), constituting an enormous up-front cost.
Thus, practical engineering control strategies for dealing with high-dimensional observational data revolves around dimensionality-reduction techniques.   Such methods, often based upon the singular value decomposition of the data, allow one to construct low-dimensional subspaces where computationally tractable controllers can be designed and implemented~\cite{Moore:1981,ERA:1985,HLBR_turb,rowley.05}.  Balanced truncation is a classic method developed to specifically take advantage of underlying low-dimensional observable and controllable subspaces to create a balanced, reduced-order model \cite{Moore:1981}.  Generalizations of this scheme which combine balanced truncation with the SVD on empirical data, such as the balanced proper orthogonal decomposition, have already been shown to overcome some of the computational difficulties associated with the high-dimension of complex systems, but still requires a pernicious linear adjoint calculation \cite{willcox:2002,rowley.05,ilak.08}.
Further innovations around system identification methods, such as the Eigensystem Realization Algorithm (ERA) and the Observer Kalman Filter Identification (OKID), were developed to aid in the discovery of input-output models for systems with control \cite{ERA:1985,OKID:1991}.  The dimension of the measurements, though, were assumed to be low and the system linear.

DMDc has a number of advantages for high-dimensional, complex systems.  First, it is based upon the DMD algorithm which is a data-driven, {\em equation-free} architecture that reconstructs the underlying dynamics of the system from snapshot measurements alone~\cite{Schmid:2009, Rowley:2009, Chen:2012,Tu:2014a}.  Substantial success has been achieved in the application of DMD to fields such as fluid dynamics which have been historically difficult to analyze and construct controllers due to the enormous number of spatial states required for simulation \cite{Grilli:2012, Schmid:2011, Schmid:2012,Bagheri:2013, Tu:2014b,tissot2014model}.  Second, DMD has acquired popularity as a method for systems with nonlinear dynamics, due to a strong connection between DMD and Koopman operator theory \cite{Koopman:1931,Mezic:2005,Rowley:2009,budivsic2012applied,Mezic:2013}.  
Finally, DMD can be modified to take advantage of sparse, or limited, measurements of the complex system~\cite{brunton:2014b,Tu:2014b,Jovanovic2014sparsity}.  Sparse measurements have recently been leveraged in a variety of complex systems, some for control~\cite{Glauser:2013,Schaeffer:2013,Lin2013ieeetac,Fardad2014ieeetac}.
Such a scenario arises in many physical, biological and engineering systems due to limited numbers of sensors.  Such advantages, in combination with the control architecture advocated here, warrant serious consideration of the DMDc as an equation-free control strategy in complex systems.

As a motivating example, DMDc can be applied to the field of computational epidemiology focusing on the eradication of diseases.  The advent of new monitoring tools and a substantial focus on the quantitative assessment of resource allocation is beginning to generate large sets of data describing the spread of infectious disease.  A substantial literature exists focused on mathematically modeling the spread of infectious disease and the effect of external control (e.g. vaccinations for Polio and bed nets for Malaria) \cite{anderson.92}.  A common challenge in computational epidemiology is deciding how to model the spread of disease leading to an enormous number of phenomenological models \cite{keeling.08}.  Equation free techniques such as DMD and DMDc provide a complementary modeling tool for analyzing the spatial-temporal spread of infectious disease.  Focusing on {\it only} the historical data containing state information (i.e. number of infections in a spatial location in a given time) and whether control interventions have been applied (i.e. number of vaccinations in a spatial location in given time), DMDc discovers the dynamical properties of the complex systems.

The outline of the paper is as follows:  \S~\ref{s:back} describes the background on the method DMD.  \S~\ref{s:DMDc} describes the new method Dynamic Mode Decomposition with control.  The following section \S~\ref{sec:applications} presents a number of numerical examples including an artificial application based on an epidemiological problem.  \S~\ref{s:systemID} discusses a number of similarities and differences from system identification methods.   

\section{Background:  Dynamic Mode Decomposition}
\label{s:back}

Dynamic Mode Decomposition (DMD) is a powerful data-driven method for analyzing complex systems.  Using measurement data from numerical simulations or laboratory experiments, DMD attempts to extract important dynamic characteristics such as unstable growth modes, resonance, and spectral properties.  This section provides the mathematical mathematical background of DMD \cite{Schmid:2009,Rowley:2009,Tu:2014a}.  

\subsection{Dynamical systems and data}

The fundamental assumption that connects the state of a linear dynamical system $\mathbf{x_{k}}$ to the next $\mathbf{x_{k+1}}$ is  
\begin{align}
\bx_{k+1} = \bA \bx_k,
\label{eq.dyn}
\end{align}
where $\bx \in \mathbb{R}^n$ and $\bA \in \mathbb{R}^{n \times n}$.  The process under observation is often continuous (whether from a numerical model or experiment), and measurements of the continuous state $\bx(t)$ can be collected at regular time intervals $\Delta t$ denoted by $\bx_k = \bx(k\Delta t)$.  Each measurement in time $\bx_k$ will be referred to as snapshots within this manuscript \cite{Sirovich:1987}.  We denote the sequence of snapshots collected by the following description: 
\begin{align}
\bX &= \left[ \begin{array}{cccc} | & |  &  & | \\
\bx_1 &\bx_{2}  & \dots & \bx_{m-1} \\
 | & |  &  & | \end{array} \right], \nonumber \\
~\bX' &= \left[ \begin{array}{ccccc} | & | &  & | \\
\bx_2 &\bx_{3}& \dots & \bx_{m~~~} \\
 | & | &  & | \end{array} \right],
\label{eq.snapshots}
\end{align}
where $m$ is total number of snapshots and $\bX'$ is the time-shifted snapshot matrix of $\bX$, i.e. $\bX' = \bA \bX$.  For DMD, data is often collected at regular time intervals $\Delta t$.  The number of snapshots required for DMD varies with the application, but is intimately related to the linearity properties of the Koopman Operator.  The solution will converge by decreasing the recording interval $\Delta t \rightarrow 0$, thus indicating the number of snapshots required (for an illuminating numerical example see \cite{Schmid:2009}).  New directions for DMD have focused on novel paradigms for collecting data in time \cite{Tu:2014a} and across the state of the system \cite{brunton:2014b}.  Each are utilizing the concepts of sparsity and compressed sensing techniques \cite{Baraniuk:2007}.   

The dynamical system Eq.~(\ref{eq.dyn}) and data snapshots Eq.~(\ref{eq.snapshots}) can be described more compactly in the following form:
\begin{align}
\bX = \bA \bX'.
\label{eq.dyn2a}
\end{align}
Solving for an approximation of the process matrix $\bA$ given the data matrices $\bX$ and $\bX'$ is the primary objective of DMD. 

\subsection{Dynamic Mode Decomposition}
\label{ss.dmdback}

The following section describes how to find the dynamic modes and eigenvalues of the underlying system $\bA$ described in Eq.~(\ref{eq.dyn2a}).  We can find $\bA$ by using the following definition:
\begin{align}
  \bA = \bX'\bX^{\dagger} ,
\label{eq.def}
\end{align}
where $^\dagger$ is the Moore-Penrose pseudoinverse \cite{Tu:2014a}.  A computationally efficient and accurate method for finding the pseudo inverse is via the singular value decomposition (SVD).  The SVD of $\bX$ results in the well-known decomposition :  
\begin{align}
\bX = \bU \bSigma \bV^* &= \left [ \tilde{\bU} ~~\tilde{\bU}_{\text{rem}} \right ] \left [ \begin{array}{cc} \tilde{\bSigma} & 0 \\ 0 & \bSigma_{\text{rem}} \end{array} \right ] \left [ \begin{array}{c} \tilde{\bV}^* \\ \tilde{\bV}_{\text{rem}}^* \end{array} \right ],\\
&\approx \tilde{\bU} \tilde{\bSigma}  \tilde{\bV}^*
\label{eq.svd}
\end{align}
where $\bU \in \mathbb{R}^{n \times n}$, $\bSigma \in \mathbb{R}^{n \times m-1}$, $\tilde{\bV}^* \in \mathbb{R}^{m-1 \times m-1}$, $\tilde{\bU} \in \mathbb{R}^{n \times r}$, $\tilde{\bSigma} \in \mathbb{R}^{r \times r}$, $\tilde{\bV}^* \in \mathbb{R}^{r \times m-1}$, $_{\text{rem}}$ indicates the remaining $m-1-r$ singular values, and $^*$ denotes the complex conjugate transpose.  Eq.~(\ref{eq.svd}) demonstrates how to reduce the dimension of the data matrix $\bX$ by appropriately choosing a truncation value $r$ of the singular values thus eliminating the remainder ($\text{rem}$) terms and allowing for the psuedo-inverse to be accomplished since $\tilde{\bSigma}$ is square.  Choosing the appropriate truncation value $r$ has a rich scientific history; notably, the Eckart-Young theorem provides a rigorous and popular method for choosing $r$ \cite{Eckart:1936,MIrsky:1960,Golub:1970}.  In addition, there are recent theoretical developments attempting to identify the correct $r$ when $\bX$ may have additive noise~\cite{Donoho:2013a,Donoho:2013b}.

Using the SVD of the snapshot matrix $\bX$ in Eq.~(\ref{eq.svd}), the following approximation of the matrix $\bA$ can be computed:
\begin{align}
\bA \approx \bar{\bA} = \bX'  \tilde{\bV} \tilde{\bSigma}^{-1}\tilde{\bU}^*,
\label{eq.dmdwoc}
\end{align}
where $\bar{\bA}$ is an approximation of the operator $\bA$ from Eq.~(\ref{eq.svd}).  A dynamic model of the process can be constructed given by the following:
\begin{align}
\bx_{k+1} =  \bar{\bA} \bx_k,
\label{eq.LargedynModel}
\end{align}
where $\bx$ and $\bar{\bA}$ have the same dimension as the matrices described earlier in Eq.~(\ref{eq.dyn}).  An eigenvalue analysis of the matrix $\bar{\bA}$ would produce the dynamic modes and eigenvalues of the system. The computation, though, can be prohibitively expensive if $n \gg 1$.   

If $r\ll n$, a more compact and computationally efficient model can be found by projecting $\bx_k$ on to a linear subspace of dimension $r$.  This basis transformation takes the form $ \mathbf{P} \bx =\tilde{\bx}$.  As previously shown by DMD, a convenient transformation has already been computed via the SVD of $\bX$, given by $\mathbf{P} = \tilde{\bU}$.  The reduced-order model can be derived as follows:  
\begin{align}
\tilde{\bx}_{k+1} &= \tilde{\bU}^* \bar{\bA}\tilde{\bU} \tilde{\bx}_k\\
 &= \tilde{\bU}^* \bX' \tilde{\bV} \tilde{\bSigma}^{-1} \tilde{\bx}_k  \\ 
 &=\tilde{\bA}\tilde{\bx}_k.
 \label{eq.dmdmodelrep}
 \end{align}
The reduced-order-model is given by the following: 
 \begin{align}
 \tilde{\bA} &= \tilde{\bU}^* \bX' \tilde{\bV} \tilde{\bSigma}^{-1}.
\label{eq.dmdwoc}
\end{align}
The eigendecomposition of $\tilde{\bA}$ defined by $\tilde{\bA}\mathbf{W} = \mathbf{W} \Lambda$ yields eigenvalues and eigenvectors that can be investigated for fundamental properties of the underlying system such as growth modes and resonance frequencies.  In addition, the computation is efficient since $\tilde{\bA} \in \mathbb{R}^{r \times r}$ and $r \ll n$.  \\
\\  
{\bf Remark} Computing the eigendecomposition of $\tilde{\bA}$ versus $\bar{\bA}$ can be a computationally crucial step for efficiency.  For example, the domain discretization of a fluids or epidemiological problem can have an arbitrarily large set of dimensions $n$.  The direct solution of the $n\times n$ eigenvalue problem might not be feasible, thus solving the $r \times r$ is substantially more attractive.  The observation is reminiscent of the Method of Snapshots by Sirovich \cite{Sirovich:1987}. 
\\
\\
For DMD, the eigenvalues of $\tilde{\bA}$ and $\bar{\bA}$ are equivalent \cite{Schmid:2009} and the eigenvectors are related via a linear transformation.  The eigenvectors of $\bar{\bA}$ are called dynamic modes \cite{Schmid:2009,Tu:2014a}.  Note, there is a difference between computing the dynamic modes with the Exact DMD method from Tu et.al. \cite{Tu:2014a} and Schmid \cite{Schmid:2009}. Here we describe the Exact DMD method giving the following relationship between the eigenvectors of $\tilde{\bA}$ and the dynamic modes $\mathbf{\phi}$ of $\bar{\bA}$:
\begin{align}
\mathbf{\phi} = \bX' \tilde{\bV} \tilde{\bSigma}^{-1}\mathbf{w}.
\label{eq.dymodes}
\end{align}

If $\lambda \neq  0$, then this is the DMD mode for $\lambda$.  If the eigenvalue is $0$, then the dynamic mode is computed using $\mathbf{\phi} = \tilde{\bU} \mathbf{w}$.  The Exact DMD algorithm has a number of advantages over the original procedure; for a detailed discussion, see \cite{Tu:2014a}.

\section{Dynamic Mode Decomposition with Control}
\label{s:DMDc}

\begin{figure*}
\begin{center}
\includegraphics[width=0.8\textwidth]{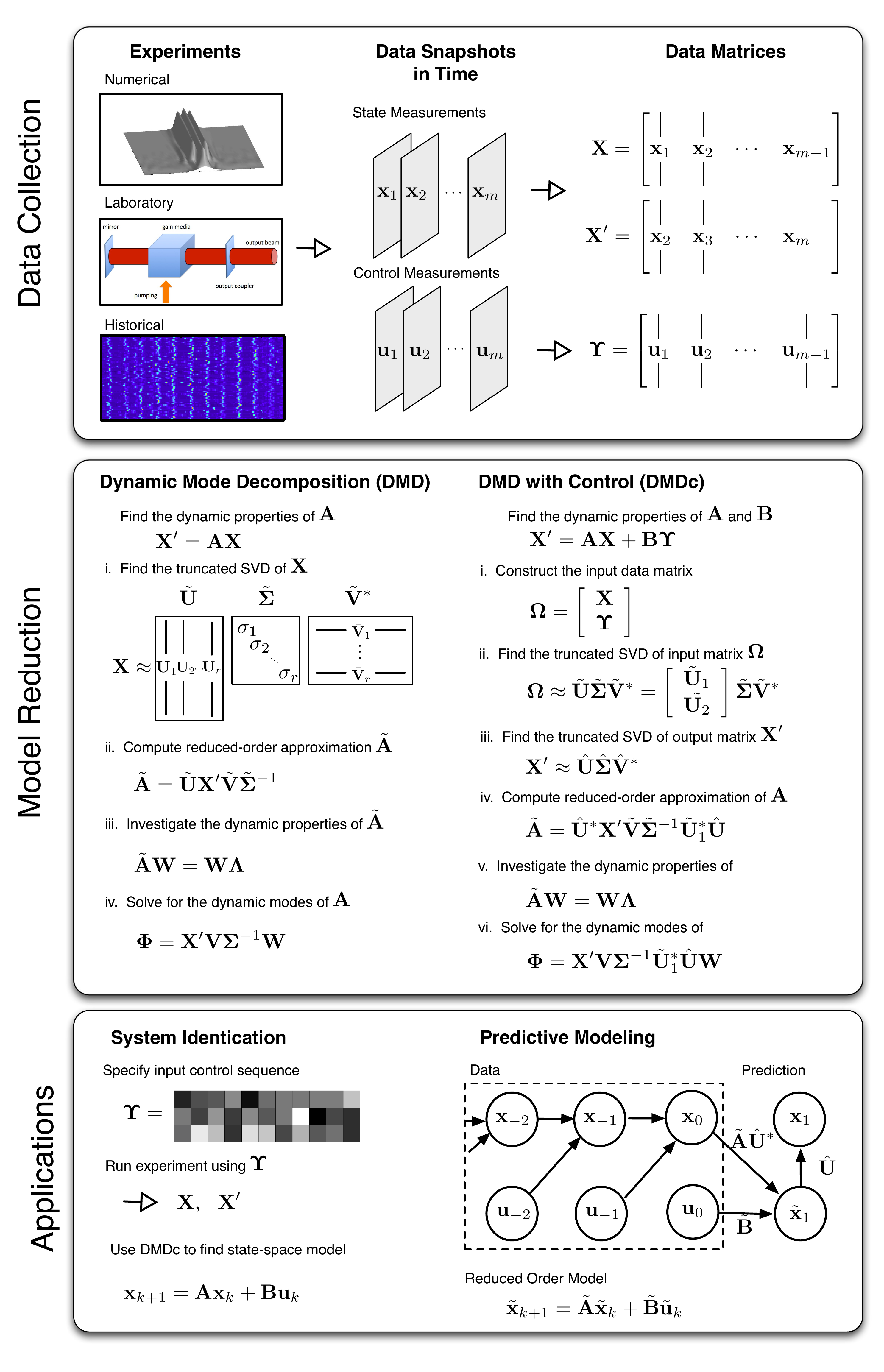}
\end{center}
\vspace{-.25in}
\caption{The illustration outlines the three major components of applying DMDc.  The top panel describes the collection of data from either a numerical, laboratory, or historical data and the curation of the data in to matrices for the methods.  Note, the figure in the historical plot is the data representing pre-vaccination Measles cases in the UK normalized similar to that found in \cite{Grenfell}.  The middle panel outlines the procedure for DMD and DMDc for comparison.  The bottom panel illustrates two practical applications of DMDc. }\label{fig:prop}
\end{figure*}

This section presents the mathematical description of Dynamic Mode Decomposition with control (DMDc).  Understanding the dynamic characteristics of complex systems that have both internal dynamics and applied external control is fundamental to controller design and sensor placement.  The DMDc method helps discover the underlying dynamics without the confounding effect of external control.  In addition, the method also quantifies the effect of control inputs on the state of the system.  Fig.~\ref{fig:prop} illustrates the data collection, the algorithm, and applications of DMDc.

The underlying dynamical system and measured data matrices are redefined to include systems with control inputs in \S~\ref{ss.dynsysc}.  The subsequent section \S~\ref{ss.known} describes how to solve for the dynamic modes if the effect of the inputs on the state is already well-known or well-estimated.  The last section \S~\ref{ss.unknown} shows how to solve for both the dynamic modes and the input matrix.

\subsection{Dynamical system with control}
\label{ss.dynsysc}
The new method modifies the basic assumption of DMD.  The linear dynamical system connecting the future state $\bx_{k+1}$ now relies on information from both the current state $\bx_{k}$ and the current control $\bu_{k}$ given by the following:  
\begin{align}
\bx_{k+1} = \bA \bx_k + \bB \bu_k,
\label{eq.dyn2}
\end{align}
where $\bx_j \in \mathbb{R}^n$, $\bu_j \in \mathbb{R}^l$, $\bA \in \mathbb{R}^{n \times n}$, and $\bB \in \mathbb{R}^{n \times l}$.  Data matrices can be constructed with temporal snapshots of the state and control input over time.  The state snapshots $\bX$ and $\bX'$ are collected in the same manner as Eq.~(\ref{eq.snapshots}).  We denote a new sequence of control input snapshots collected by the following description: 
\begin{align}
\bUp &= \left[ \begin{array}{cccc} | & |  &  & | \\
\bu_1 &\bu_{2}  & \dots & \bu_{m-1} \\
 | & |  &  & | \end{array} \right].
\label{eq.snapshots2}
\end{align}
Eq.~(\ref{eq.dyn2}) can be rewritten to include the new data matrices:  
\begin{align}
\bX' = \bA \bX + \bB \bUp.
\label{eq.dyn3b}
\end{align}
Utilizing the three data matrices, approximations of the linear mappings $\bA$ and $\bB$ can be found.  In the following two sections, we describe how to find the dynamic modes of $\bA$ given the inclusion of control snapshots.  The first section outlines the analysis and algorithm if the matrix $\bB$ is known or well estimated.  If unknown, the second section describes how to discover both $\bA$ and $\bB$ from the observation matrices.

\subsection{The map $\bB$ is known}
\label{ss.known}

The following section describes how to find the dynamic modes and eigenvalues of the underlying system $\bA$ when the matrix $\bB$ is known.  The assumption that $\bB$ is known or well-estimated is an idealistic view of most complex systems, but it helps provide one of the major motivations for this work.  Finding the underlying dynamics $\bA$ in a complex system where control has been applied is essential for designing controllers and placement of sensors.  If external control has been applied to the system, standard DMD would produce incorrect dynamic information.  The more general case where $\bB$ is unknown will be described in the following section.

Eq.~(\ref{eq.dyn3b}) can be re-arranged by pairing the time-shifted state snapshot matrix with the control snapshot matrix and the known matrix $\bB$.  
\begin{align}
\bX' - \bB \bUp = \bA \bX 
\label{eq.dyn3}
\end{align}
The mapping $\bA$ can be solved for similar to Eq.~(\ref{eq.def}).  Again, the truncated singular value decomposition of $\bX$ gives the matrix factorization $\tilde{\bU}\tilde{\bSigma}\tilde{\bV}^*$.  Thus, the approximation of $\bA$ is given by the following description:  
\begin{align}
\bA \approx \bar{\bA} = ( \bX'- \bB \bUp) \tilde{\bV} \tilde{\bSigma}^{-1}\tilde{\bU}^*.
\label{eq.dmdwc}
\end{align}
Note, if the control snapshots are $\bu_j = \mathbf{0},~\forall~ j~ \in [1,m]$, then the derivation is equivalent to DMD.  A dynamic model of both the computed process and the given input matrix can be constructed described by the following:
\begin{align}
\bx_{k+1} =  \bar{\bA} \bx_k+ \bB \bu_k 
\label{eq.LargedynModel}
\end{align}
where $\bx$, $\bar{\bA}$, and $\bB$ are the same dimensions of the matrices described earlier in Eq.~(\ref{eq.dyn}).  If $r\ll n$ though, a more compact and computationally efficient model can be found using the same basis transformation $\mathbf{P} \bx = \tilde{\bx}$ as described earlier for DMD.  Again, a convenient transformation has already been computed via the SVD of $\bX$, given by $\mathbf{P} = \tilde{\bU}$.  The reduced-order model can be derived as follows:  
\begin{align}
\tilde{\bx}_{k+1} &= \tilde{\bU}^* \bar{\bA}\tilde{\bU} \tilde{\bx}_k+ \tilde{\bU}^* \bB \bu_k,\\
 &= \tilde{\bU}^* (\bX' - \bB\bUp )\tilde{\bV} \tilde{\bSigma}^{-1} \tilde{\bx}_k + \tilde{\bU}^*\bB \bu_k \\ 
 &=\tilde{\bA}\tilde{\bx}_k+\tilde{\bB} \bu_k
 \label{eq.dmdmodelrep}
 \end{align} 
The reduced-order approximation of $\bA$ is given by the following:
 \begin{align}
 \tilde{\bA} &= \tilde{\bU}^*(\bX' - \bB\bUp ) \tilde{\bV} \tilde{\bSigma}^{-1} 
\label{eq.dmdwoc}
\end{align}
The eigendecomposition of $\tilde{\bA}$ defined by $\tilde{\bA}\mathbf{W} = \mathbf{W}\mathbf{\Lambda}$ yields eigenvectors that can be used to find the dynamic modes.  Similar to Exact DMD, the dynamic modes can be found with the following description: 
\begin{align}
\mathbf{\phi} = ( \bX' - \bB \bUp) \tilde{\bV} \tilde{\bSigma}^{-1}\mathbf{w}
\label{eq.dymodes}
\end{align}
If $\lambda \neq  0$, then this is the DMD mode for $\lambda$.  If the eigenvalue is $0$, then the dynamic mode is computed using $\mathbf{\phi} = \tilde{\bU}\mathbf{w}$.
 
\subsection{The map B is unknown}
\label{ss.unknown}
The assumption that $\bB$ is known indicates a significant amount of knowledge about how control affects the system.  This section relaxes that assumption and notably demonstrates that approximations of the matrices $\bA$ and $\bB$ can both be found from state and control snapshots.  To the experimentalist or analyst, this is by far more interesting since only the snapshots of the control and state are required to find the properties of the underlying process $\bA$ and how that process is affected by control $\bB$.  

The dynamical system from Eq.~(\ref{eq.dyn3b}) can be manipulated giving the following representation:
\begin{align}
\bX' = \left[\bA ~~\bB \right] \left[ \begin{array}{c} \bX \\ 
 \bUp \end{array}\right] = \mathbf{G}  \mathbf{\Omega},
\label{eq.noB2}
\end{align}
where $\mathbf{\Omega}$ contains both the state and control snapshot information.  Here, we again seek a best-fit solution of the operator $\mathbf{G}$ which now contains the process dynamics $\bA$ and input matrix $\bB$.  A SVD is performed on the augmented data matrix giving $\mathbf{\Omega} =\bU\bSigma \bV^* \approx \tilde{\bU} \tilde{\bSigma} \tilde{\bV}^* $.  The truncation value of the SVD for $\mathbf{\Omega}$ will be defined as $p$.  Note, the truncation value of $\mathbf{\Omega}$ should be larger than of $\bX$.  The following computation provides an approximation of $\mathbf{G}$:
\begin{align}
\mathbf{G} \approx \mathbf{\bar{G}} =  \bX'  \tilde{\bV} \tilde{\bSigma}^{-1}\tilde{\bU}^*,
\label{eq.dmdcF}
\end{align}
where $\mathbf{G} \in \mathbb{R}^{n \times (n+l)}$.  We can now find approximations of the matrices $\bA$ and $\bB$ by breaking the linear operator $\tilde{\bU}$ in to two separate components given by the following:
\begin{align}
[\bA,~~\bB] &\approx [\bar{\bA},~~\bar{\bB}] \\
&\approx [\bX' \tilde{\bV} \tilde{\bSigma}^{-1} \tilde{\bU}_1 ^*,~~\bX' \tilde{\bV} \tilde{\bSigma}^{-1} \tilde{\bU}_2 ^*]
\label{eq.dmdcF2}
\end{align}
where $\tilde{\bU}_1 \in \mathbb{R}^{n \times p}$, $\tilde{\bU}_2 \in \mathbb{R}^{l \times p}$, and $\tilde{\bU} = [\tilde{\bU}_1^*~~ \tilde{\bU}_2^*]^T$. Similar to Eq.~(\ref{eq.dmdcF2}), a dynamic model using the matrices $\bar{\bA}$ and $\bar{\bB}$, but for a large dimensional system where $n\gg1$, this is computationally prohibitive. Here, we again seek a reduced order model of rank $r\ll n$ where a transformation is required such that $\bx = \mathbf{P}\tilde{\bx}$ and $\tilde{\bx} \in \mathbb{R}^{r}$.  

Unlike DMD, the truncated left singular vectors $\tilde{\bU}$ can not be used to define the subspace on which the state evolves.  For Eq.~(\ref{eq.dmdcF2}), the truncated left singular vectors of $\mathbf{\Omega}$ define the {\bf input} space.  To find a linear transformation $\mathbf{P}$ for the state $\bx$, we utilize a reduced-order subspace of the {\bf output} subspace.  This fundamental observation allows for DMDc to discover a reduced-order representation of the dynamics $\bA$ and input matrix $\bB$. 

To find the reduced-order subspace of the output space, a second singular value decomposition is required.  The data matrix of the {\bf output} space $\bX'$ can be approximated by the familiar SVD: $\hat{\bU} \hat{\bSigma} \hat{\bV}^*$ where the truncation value is $r$ and $\hat{\bU} \in \mathbb{R}^{n \times r}$, $\hat{\bSigma} \in \mathbb{R}^{r \times r}$, and $\hat{\bV}^* \in \mathbb{R}^{r \times {m-1}}$.  Note, the two SVDs will likely have different truncation values of the input and output matrices $p$ and $r$ and $p>r$.  Using the transformation $\bx = \hat{\bU} \tilde{\bx}$, the following reduced-order approximations of $\bA$ and $\bB$ can be computed:

\begin{align}
\tilde{\bA} &= \hat{\bU}^* \bar{\bA} \hat{\bU} =  \hat{\bU}^* \bX' \tilde{\bV} \tilde{\bSigma}^{-1} \tilde{\bU}_1 ^* \hat{\bU} \\
\tilde{\bB} &= \hat{\bU}^* \bar{\bB} ~~~= \hat{\bU}^* \bX' \tilde{\bV} \tilde{\bSigma}^{-1}\tilde{\bU}_2 ^* 
\label{eq.dmdcF3}
\end{align}
where $\tilde{\bA} \in \mathbb{R}^{r \times r}$ and $\tilde{\bB} \in \mathbb{R}^{r \times l}$.  We can then form the reduced order equation as Eq.~(\ref{eq.dmdmodelrep}) given by the following
\begin{align}
\tilde{\bx}_{k+1} = \tilde{\bA}\tilde{\bx}_k +\tilde{\bB} \bu_k 
\label{eq.dmdcF3}
\end{align}
Similar to DMD, the dynamic modes of $\bA$ can be found by first solving the eigenvalue decomposition $\tilde{\bA} \mathbf{W} = \mathbf{W} \mathbf{\Lambda}$.  The transformation from eigenvectors to dynamic modes of $\bA$ is slightly modified and is given by the following:   
\begin{align}
\mathbf{\phi} =  \bX' \tilde{\bV} \tilde{\bSigma}^{-1} \tilde{\bU}_1 ^*\hat{\bU} \mathbf{w}.
\label{eq.dymodes2}
\end{align}
where the relationship between $\mathbf{\phi}$ and $\mathbf{w}$ is similar to Exact DMD.  

\subsection{The algorithm}
\label{ss.algo}
The following section outlines the algorithm.
\begin{enumerate}
\item {\bf Collect and construct the snapshot matrices:}
\\
Collect the state and control snapshots and form the matrices $\bX$, $\bX'$, and $\bUp$ as described in Eq.~(\ref{eq.snapshots}) and Eq.~(\ref{eq.snapshots2}).  Stack the data matrices $\bX$ and $\bUp$ to construct the matrix $\mathbf{\Omega}$.

\item {\bf Compute the SVD of the input space $\mathbf{\Omega}$.}\\
Compute the singular value decomposition of $\mathbf{\Omega}$ as described in Eq.~(\ref{eq.svd}) thereby obtaining the decomposition $\mathbf{\Omega} \approx \tilde{\bU} \tilde{\bSigma} \tilde{\bV}^*$ with truncation value $p$.

\item {\bf Compute the SVD of the output space $\bX'$.}\\
Compute the singular value decomposition of $\bX'$ as described in Eq.~(\ref{eq.svd}) thereby obtaining the decomposition $\bX' \approx \hat{\bU} \hat{\bSigma} \hat{\bV}^*$ with truncation value $r$.

\item {\bf Compute the approximation of the operators $\mathbf{G}=[\bA~~\bB] $}
\\
Compute the following:
\begin{align}
\tilde{\bA} &= \hat{\bU}^* \bX' \tilde{\bV} \tilde{\bSigma}^{-1} \tilde{\bU}_1^* \hat{\bU}\\
\tilde{\bB} &= \hat{\bU}^* \bX'  \tilde{\bV} \tilde{\bSigma}^{-1} \tilde{\bU}_2^*
\end{align}
\item {\bf Perform the eigenvalue decomposition of $\tilde{\bA}$}
\\
Perform the eigenvalue decomposition given by the following:
\begin{align}
\tilde{\bA}\mathbf{W} = \mathbf{W}\mathbf{\Lambda}
\end{align}

\item {\bf Compute the dynamic Modes of the operator A}
\\
\begin{align}
\mathbf{\Phi} =  \bX' \tilde{\bV} \tilde{\bSigma}^{-1} \tilde{\bU}_1^*\hat{\bU}  \mathbf{W}
\label{eq.dymodesF}
\end{align}
\end{enumerate}
\section{Applications}\label{sec:applications}

The following section describes a number of numerical examples for the application of this method.  The examples increase in complexity as the section progresses.  The emphasis for each of these examples is the benefit of including control snapshot information to the analysis.  

\subsection{Example 1 -- Unstable linear system with proportional controller}\label{sec:results:unlinearsys}

DMDc can help discover the underlying dynamics of a system through measurements of both the state and external inputs.  Here, we demonstrate the idea on a simple two-dimensional unstable linear system with a stabilizing controller.  Despite the simplicity of the mathematical problem, the example is illustrative for the general concept of DMDc.   Consider the following dynamical system: 
\begin{align}
\left[ \begin{array}{c} x_1 \\ x_2 \end{array} \right ]_{k+1} = \left[ \begin{array}{cc} 1.5 & 0 \\0 & 0.1 \end{array}\right]  \left[ \begin{array}{c} x_1 \\ x_2 \end{array} \right ]_{k} + \left[ \begin{array}{c} 1 \\ 0 \end{array} \right ] u_k
\label{eq.numex1}
\end{align}
where $u_k = K [x_1]_k$ and $K = -1$.  The proportional controller clearly stabilizes the system by moving the unstable eigenvalue within the unit circle.  If we have access to the input data and the $\bB$ matrix as described in \S\ref{ss.unknown}, we can collect state and control snapshots to perform the DMDc computation.  For an initial condition $[4~~7]^T$, the following are the data matrices constructed from computing the first five temporal snapshots of Eq.~(\ref{eq.numex1}):
\begin{align}
\bX &= \left[ \begin{array}{cccc} 4 & 2 & 1 & 0.5  \\ 7 & 0.7 & 0.07 & 0.007  \end{array} \right] \\
\bX' &= \left[ \begin{array}{cccc} 2 & 1 & 0.5 & 0.25 \\  0.7 & 0.07 & 0.007 & 0.0007 \end{array} \right] \\
\bUp &= \left [\begin{array}{cccc} -4 &-2 &-1 &-0.5 \end{array} \right ]
\label{eq.numexsnap}
\end{align}
Following the description in \S\ref{ss.known}, we compute the singular value decomposition of $\bX$.  Here, we use MATLAB's economy sized singular value decomposition algorithm to give the following matrix factorization of $\bX$.  
\begin{align}
\tilde{\bU} &= \left[ \begin{array}{cc} -0.5239 & -0.8462 \\ -0.8462 & 0.5329 \end{array} \right] \\
\tilde{\bSigma} &= \left [\begin{array}{cc}8.2495 & 0 \\ 0 &1.6402  \end{array} \right ] \\
\tilde{\bV} &=\left [\begin{array}{cc}-0.9764 & 0.2105 \\ -0.2010 &-0.8044 \\-0.0718 & -0.4932 \\ -0.0330 & -0.2557 \end{array} \right ]
\label{eq.numexsnap2}
\end{align}
Now, we can compute Eq.~(\ref{eq.dmdwc}) using the data matrices in Eq.~(\ref{eq.numexsnap}), the SVD matrices in Eq.~(\ref{eq.numexsnap2}), and the matrix $\bB$ in Eq.~(\ref{eq.numex1}) giving the following approximation to $\bA$:
\begin{align}
\bar{\bA}= \left[ \begin{array}{cc} 1.5 & 0 \\ 0 & 0.1\end{array} \right] 
\label{eq.approxA}
\end{align}
where we recover the unstable linear dynamics from data of the state and control snapshots.  This example demonstrates the utility of DMDc with recovering unstable dynamics from a system that would otherwise appear to be stable.

\subsection{Example 2 -- Large-scale, stable Linear Systems}\label{sec:results:linearsys}

In this section, we investigate stable linear systems where the number of measurements are significantly greater than the dimensionality of the underlying system.  The previous example demonstrated the utility of the method on a low-dimensional unstable model.  Here, the method is applied to large-scale dynamical systems that have an underlying low-dimensional attractor.  

To construct these large-scale systems, a low-dimensional stable model is generated and subsequently embedded in to a higher dimensional subspace.   There are three steps for generating the model and data matrices to compare the output of DMDc and the generated model: \\
\begin{enumerate}
  \item {\bf Generate a low-dimensional stable state-space model, $\bA$ and  $\bB$} \\
Generate discrete random state-space systems using MATLAB's command Discrete Random State Space Method.  These stable-discrete state-space models can be used as numerical experiments for DMDc. Here, we have chosen a 5 dimensional model, 2 input variables, and 100 measurement variables.  The output is a state space model $\tilde{\bA}$, $\tilde{\bB}$, and $\mathbf{C}$.  
  \item {\bf Generate random input data $\bUp$ }\\
  Using MATLAB's randn command, generate a matrix of random inputs, $\bUp\in \mathbb{R}^{2 \times m-1}$.
  \item {\bf Use the model and input vector to generate the data matrices $\bX$ and $\bX'$} \\
  Using the model and the input matrix, generate output data for the snapshot matrix.
\end{enumerate}

Using the data matrices $\bX$, $\bX'$ and $\bUp$, the DMDc computation can be performed to find an approximation of $\tilde{\bA}$ and $\tilde{\bB}$.  To compare the generated model and the model produced by DMDc, we assign $\tilde{\mathbf{C}} = \hat{\bU}$.  The assignment allows for the comparison of state-space models.  

\begin{figure}
\begin{center}
\includegraphics[width=0.45\textwidth]{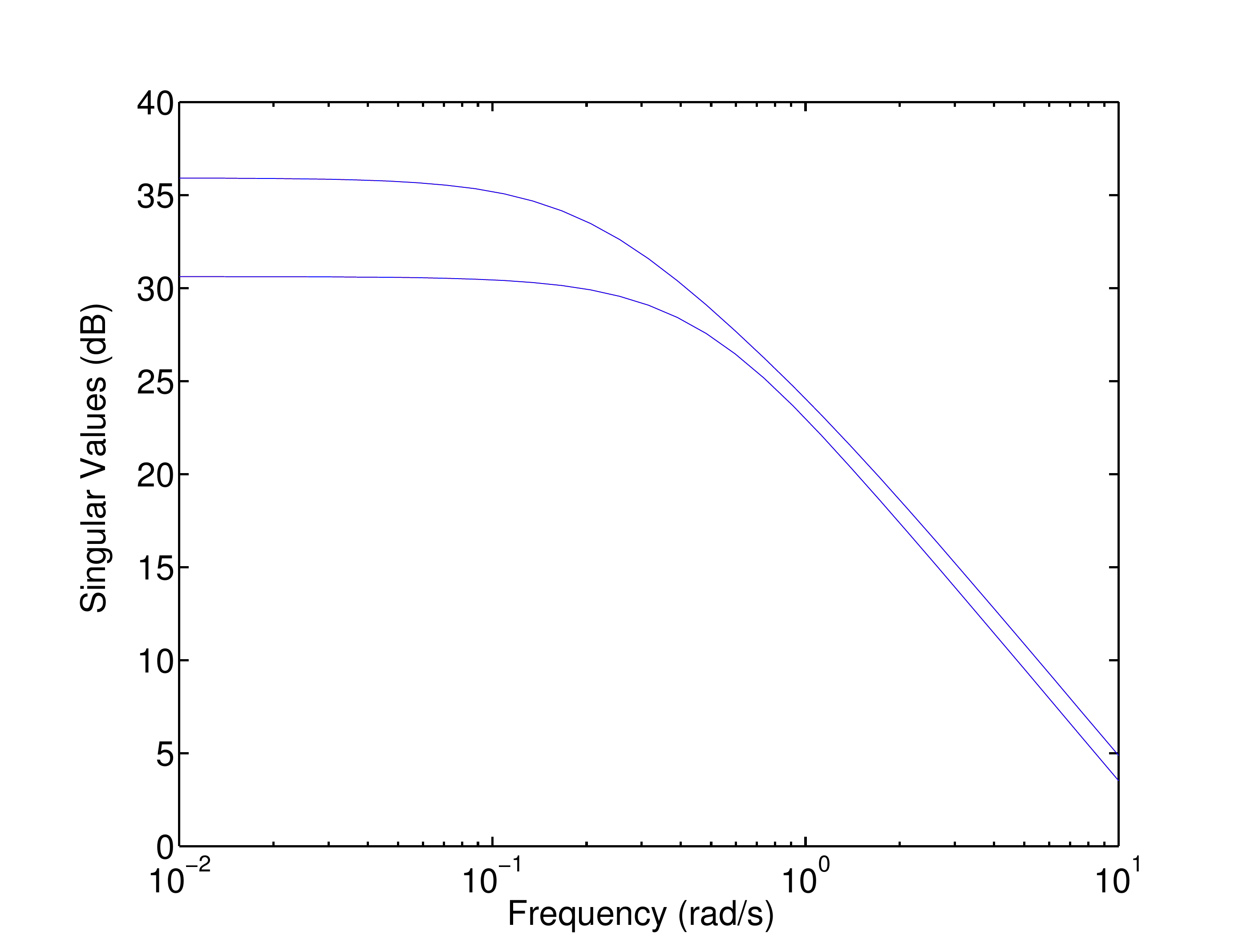}
\end{center}
\caption{The singular values of the frequency response for a large scale, stable linear systems.  The blue line is from the model from DMDc and the red is from the real model. Note, an equivalent frequency response can not be constructed from DMD alone since it does not consider input-output systems. }\label{fig:sigma}
\end{figure}

The singular values of the frequency response, a multi-input multi-output (MIMO) generalization of a BODE plot, is used to compare the two models.  The MATLAB command sigma will generate the frequency response for both systems.  Fig.~\ref{fig:sigma} illustrates one such comparison arising from a single numerical realization from the ensemble. Note, there is no distinction between the generated model (in red) and the model from DMDc (in blue) for both control inputs (both lines).   

\subsection{Example 3 --  A sparse linear system in the Fourier domain}\label{sec:results:linearfourier}
The final example for DMDc is a large-scale dynamical system on a spatial grid.  The system consists of high-dimensional full-state measurements, although the dynamics are governed by a low-dimensional dynamical system in the Fourier domain.  The motivation for this example comes from epidemiology and infectious disease spread where the measurements can be high-dimensional in both space and time.  For example, consider the number possible states of a dynamical system to represent flu infections across the world over a decade, including both spatial discretization and disease heterogeneity factors.  In this example, the underlying attractor could be quite low dimensional.  To complicate this picture, actuation in the form of a spatial delivery of vaccinations is also occurring each year, which can directly affect the dynamics of an infectious disease.  

Here, we construct a sparse dynamical system in a two-dimensional Fourier domain as an abstraction of the problem described above.  Only 5 modes are allowed to be non-zero.  The dynamical system on these spatial modes is constructed in the following way:  for each mode, a temporal oscillation frequency is chosen randomly and a small, stable damping rate is similarly chosen.  The boundary conditions are periodic, thus restricting the dynamics to a torus.  This system was previously constructed in \cite{brunton:2014b} to demonstrate compressive DMD.  Here though, the example is extended to allow for actuation in the spatial domain.  The spatial actuation is then Fourier transformed in order to compute the effect on the underlying dynamical system.  The spatial grid used is $128 \times 128$.  

\begin{figure*}
\begin{center}
\includegraphics[width=1.05\textwidth]{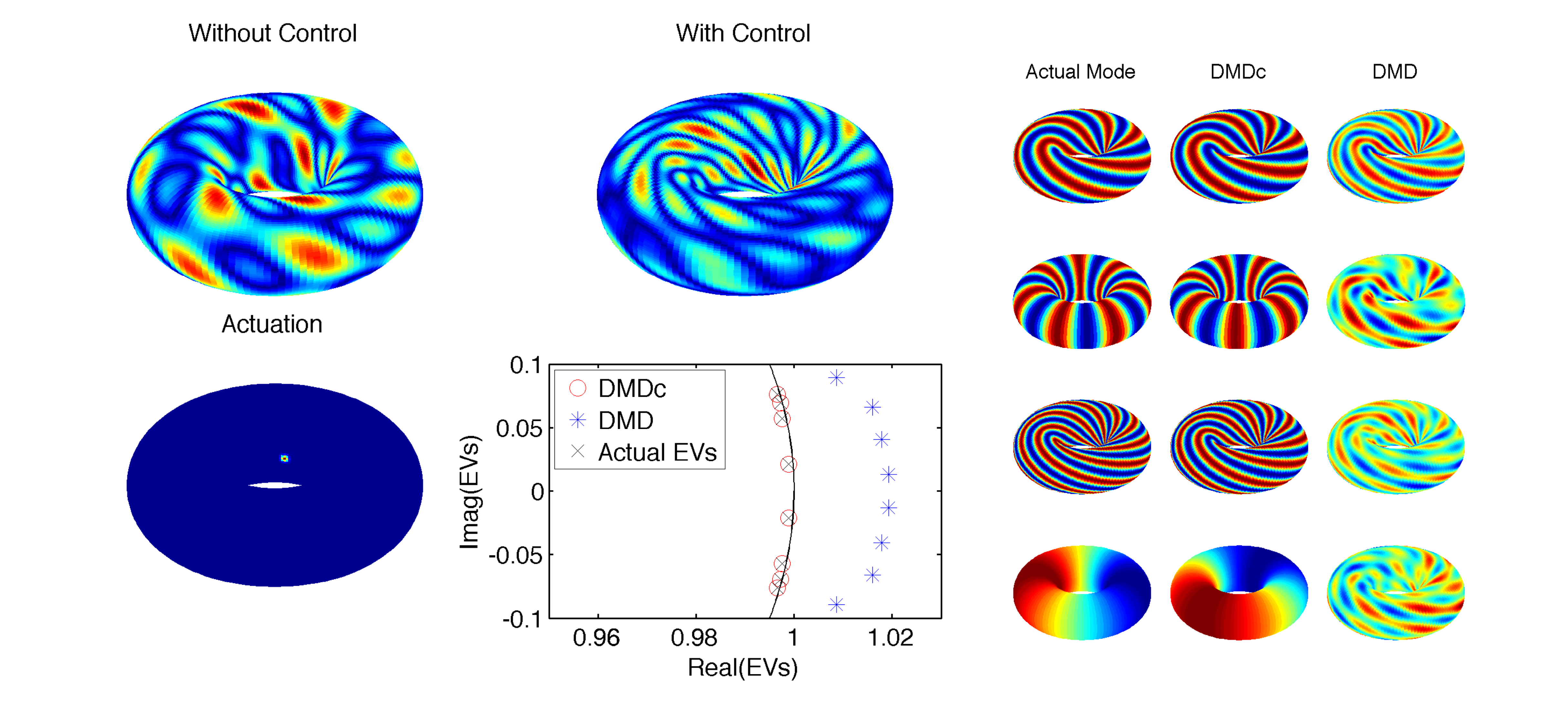}
\end{center}
\vspace{-.25in}
\caption{The top left panel illustrates one realization of Example 3 without actuation over time. The right panel illustrates the same dynamical system, but with actuation.  The bottom left panel illustrates the actuation applied in the spatial domain.  The bottom middle panel shows a comparison between the actual eigenvalues and the eigenvalues found from DMD and DMDc. On the right the first four dynamic modes of DMD and DMDc are compared to the actual underlying spatial modes. }\label{fig:finalexample}
\end{figure*}

Similar to the previous examples, the underlying dynamics of the system can be discovered soley from state and control snapshots in the spatial domain using DMDc.  The top left plot of Fig.~\ref{fig:finalexample} shows the evolution of one such unforced system in space.  The right plot shows the effect of actuation on the same system.  The actuation is a localized negative control input applied in the spatial domain, shown in the lower left plot.  The eigenvalue plot shows that DMDc discovers the underlying eigenvalues more accurately than DMD.  In addition, the zero-valued Fourier modes can be contaminated with Gaussian noise without a qualitatively change in  the behavior of DMDc.
  
\section{Connections to system identification methods}
\label{s:systemID}

This section explores the connection of DMDc to two system identification methods:  the Eigensystem Realization Algorithm (ERA) and the Observer Kalman Filter Identification (OKID).  These system identification methods were developed to derive a state-space model for control in aerospace applications involving flexible structures~\cite{ERA:1985,OKID:1991}.  The identification process involves applying control and observing system behavior.  DMDc and other modal decomposition methods can be used similarly, but may also be applied to historical data records from many other fields such as epidemiological modeling.  Here, we briefly describe the similarities and differences between the methods. 

System identification methods such as ERA/OKID were developed for input-output systems which typically have a higher rank/dimensionality than the number of observables $r > n$ \cite{kalman:1965, ERA:1985}.  To contrast, modal decomposition methods such as DMD, DMDc, Proper Orthogonal Decomposition (POD), and Balanced Proper Orthogonal Decomposition (BPOD) are typically applied to complex systems where the number of measurements are significantly larger than the rank of the underlying attractor $n \gg r$, e.g. fluid dynamics problems.  Fig.~\ref{fig:dim} illustrates the regime of applications where each of these methods are typically applied. In addition, DMD and POD have been previously established as analysis methods for {\it nonlinear} complex systems \cite{Koopman:1931,Rowley:2009,Mezic:2013}

Previous work by Tu et al. \cite{Tu:2014a} has established a number of connections between DMD and ERA.  The similarities and differences between DMDc and ERA listed in this section, though, are more readily compared since both algorithms assume input-output systems.  The following list offers a brief comparison between how DMDc and ERA differ for the construction of a typical input-output model:
 \begin{figure}
\begin{center}
\includegraphics[width=0.45\textwidth]{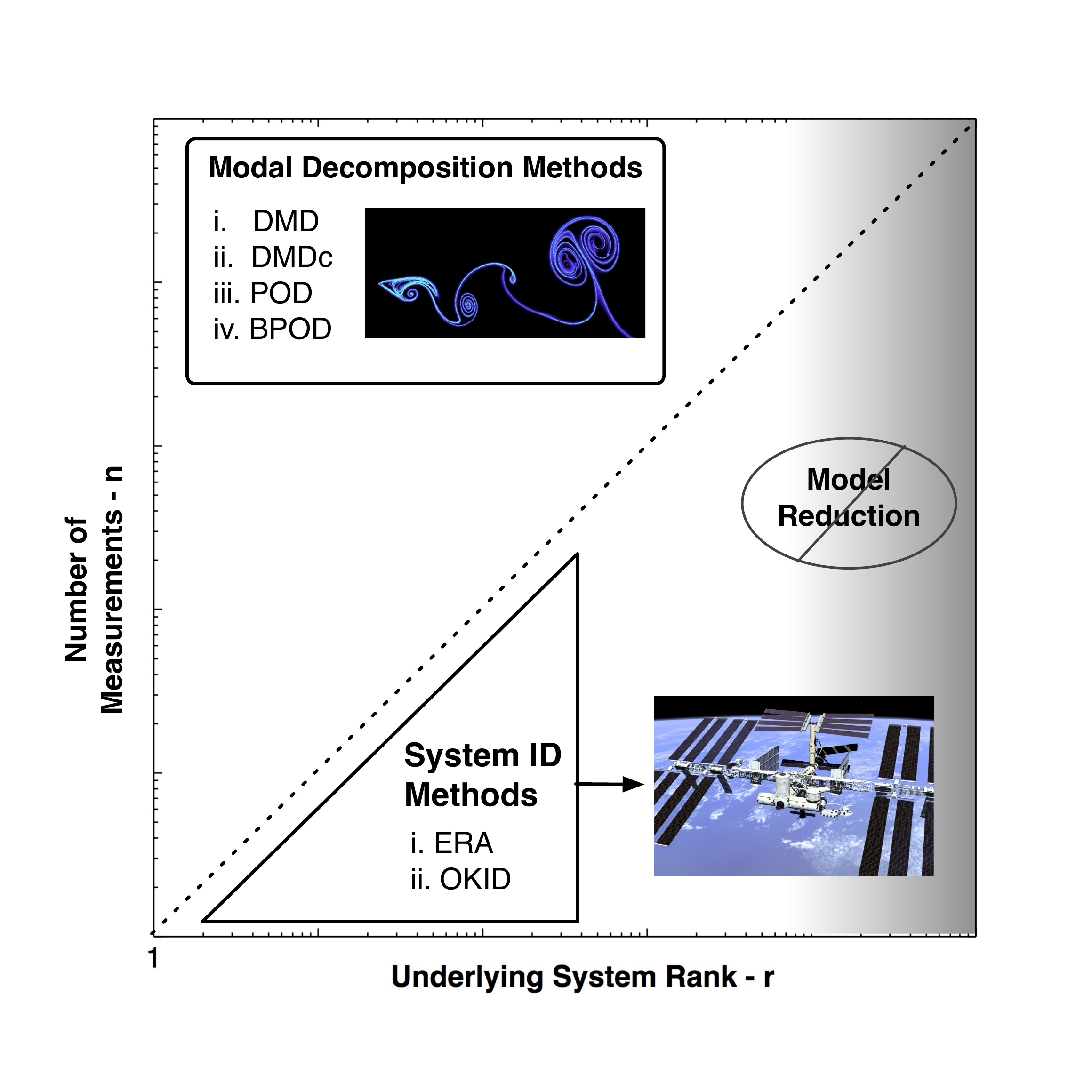}
\end{center}
\vspace{-.25in}
\caption{An illustration depicting the different regimes, with respect to the rank of the system and the number of measurements, of the modal decomposition methods and the system identification methods.}\label{fig:dim}
\end{figure}
\begin{itemize}
\item{{\bf The data matrix construction:}\\   The data from the ERA procedure is fundamentally impulse-response data, whereas DMDc can have arbitrary input histories.  The input histories are fundamental to the DMDc procedure.  Despite this difference, the construction of the state data matrices are similar between DMDc and ERA.  Both the matrix $\bX$ from DMDc and the Hankel matrix $\mathbf{H}$ of ERA assume snapshots of the state at regular intervals.  The data matrix $\mathbf{H}$ is also vertically stacked with time shifted versions of the snapshot.  DMDc does not require shift-stacking the matrix since there is little risk of column-rank deficiency due to the typically large number of observables.  The two state matrices are equivalent on the condition that the $\mathbf{H}$ is not vertically stacked with state snapshots \cite{Tu:2014a}.}

\item{{\bf The $\bA$ matrix:}\\It was previously shown that if the data matrices described in the first bullet are the same, the matrices $\bA$ produced by DMD and ERA are equivalent up to a similarity transformation \cite{Tu:2014a}.  The $\bA$ matrix constructed by DMDc is different from DMD and ERA since the {\it input} space and thus the subspace $\tilde{\bU_1^*}$ contains added information from the control snapshots.}

\item{{\bf The $\bB$ matrix }\\  To compute the matrix $\bB$ using ERA, only the first data snapshot after the impulse is utilized, which translates to the following discrete dynamical system relationship $\mathbf{x_1}=\mathbf{C}\bB\mathbf{u_0}$ where $\mathbf{u_0}$ is an impulse and $\mathbf{C}$ is the standard linear map for the observable equation.  Note, ERA also requires a projection of the single data snapshot on to the left singular vectors of the data matrix $\mathbf{H}$ and the same similarity transform $\bSigma^{1/2}$ described for the matrix $\bA$. The ERA formulation can be contrasted with DMDc through the illustration of the difference in the data matrix construction for DMDc given by the following:
\begin{align}
\mathbf{\Omega}
 = \left[ \begin{array}{ccccc}| &|&|&&| \\
  \mathbf{0}& \mathbf{x_1}&\mathbf{x_2}& \hdots& \mathbf{x_{m_c}} \\
 | &|&|&&| \\ 
\mathbf{u_0} & \mathbf{0} & \mathbf{0} & \hdots& \mathbf{0} \end{array} \right ]
\label{eq.Vaug}
\end{align}
The computation to find the matrix $\bB$ is quite different for DMDc.  Arbitrary control histories can be included in $\mathbf{\Omega}$ to compute $\bB$ whereas ERA is primarily impulse response focused.  Further, finding the matrix $\bB$ with ERA will not be as robust to noise compared with using DMDc and a longer input history.
}
\item{{\bf The $\mathbf{C}$ matrix } For DMDc, ERA, and DMD, a linear transformation matrix maps the model state to the observables.   Each of these methods utilize the left singular vectors of their data matrices for the mapping. There is an important distinction between the role of the left singular vectors for DMDc and ERA.  The mapping for DMDc projects a high-dimensional set of observables  on to a lower-dimensional subspace.  In ERA, the left singular vectors often lifts the dimension of the observables, see Fig.~\ref{fig:dim} for an illustration of the rank of the model versus dimension of the observables.  }
\end{itemize}

The observer/Kalman filter identification method allows minimal realization algorithms such as ERA to be generalized from impulse response data to data that is driven by rich input signals~\cite{OKID:1991}.  The calculation of the above matrices $\bA$, $\bB$, and $\mathbf{C}$ is typically considered more robust when combining OKID with ERA.  An often cited computational challenge confronting ERA is the analysis of lightly damped systems.  The magnitude of data (number of snapshots) may be prohibitively large for lightly damped systems, and factoring the Hankel matrix using the SVD is computationally prohibitive.  A major similarity between OKID and DMDc is the construction of the data matrix; OKID constructs an augmented data matrix that also stacks the control with the state.  Similar to DMD and ERA, in the limit of only evaluating the first row of the augmented Hankel matrix, the data matrices between DMDc and OKID are equivalent.  

\section{Discussion}\label{sec:discussion}

Complex, high-dimensional data has become ubiquitous in traditional scientific and engineering applications as well as modern data-rich fields such as internet traffic, distribution systems, and transportation networks.  Machine-learning and statistical methods have been successfully applied to characterize many of these so-called {\it big-data} problems.  Similarly, scientific and engineering fields, exemplified by control theoretic community, have focused on the development of quantitative and automatic dimensionality reduction methods to both characterize and {\it control} complex systems.  In order to construct effective controllers, the underlying system needs to be well-understood.  Accurately describing the underlying system is a challenge when the system is complex, high-dimensional, and without well-characterized governing equations.

Dynamic Mode Decomposition (DMD) is a data-driven, equation-free method that helps meet a number of these modern-day challenges.  The method has strong connection to nonlinear operator theory and discovers spatial-temporal coherent modes from data.  DMD, though, does not produce accurate reduced-order-models from complex systems with exogenous forcing.  Dynamic mode decomposition with control (DMDc) inherits the advantages of DMD, but also provides accurate input-output models for complex systems with actuation.  The method can be applied to data from a variety of sources including historical, experimental, and black-box simulations.  

Methods such as DMDc will play an increasing role in the analysis of large-scale datasets from complex systems.  DMD has already been applied to a significant number of applications in the fluid dynamics community \cite{ Schmid:2011, Schmid:2012,Grilli:2012,Bagheri:2013, Tu:2014b,tissot2014model} and is expanding to a variety of other applications like background subtraction in video processing \cite{Grosek2014}.  We believe DMDc is poised to similarly excel as a tool for a diverse set of engineering applied science applications where control of the complex system is important.  Further, the DMDc method is well-suited to couple with innovative sparsity-promoting sampling and control strategies \cite{Candes:2006,Schaeffer:2013,Glauser:2013,Fardad2014ieeetac}.  This connection has already been demonstrated for DMD both in time and space \cite{brunton:2014b,Jovanovic2014sparsity,Tu:2014b}.  DMDc is therefore positioned to have a dramatic effect on the analysis and control of large-scale complex systems.      

\section*{Acknowledgements} 

The authors would like to thank Bill and Melinda Gates for their active support of the Institute for Disease Modeling and their sponsorship through the Global Good Fund.  J. N. Kutz acknowledges support from the U.S. Air Force Office of Scientific Research (FA9550-09-0174).  Productive discussions about Dynamic Mode Decomposition with Bing Brunton, Clancy Rowley, and Jonathan Tu are likewise greatly appreciated.

\bibliographystyle{plain}
\bibliography{DMDc}
\newpage

\begin{wrapfigure}{l}{0.15\textwidth}
\vspace{-.225in}
  \begin{center}
    \includegraphics[width=.15\textwidth]{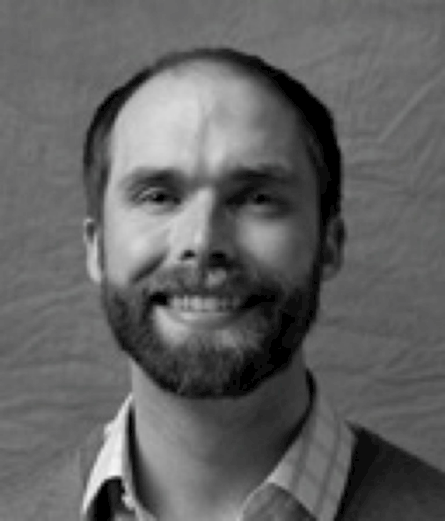}
  \end{center}
  \vspace{-.275in}
\end{wrapfigure}
\noindent\textbf{Joshua L. Proctor} received the B.S. degrees in aeronautics and astronautics, the B.A. degree in english literature from the University of Washington, Seattle, WA, in 2006, and the Ph.D. degree in mechanical and aerospace engineering from Princeton University, Princeton, NJ, in 2011.  He is currently an Associate Principal Investigator at the Institute for Disease Modeling. \\

\begin{wrapfigure}{l}{0.15\textwidth}
\vspace{-.225in}
  \begin{center}
    \includegraphics[width=.15\textwidth]{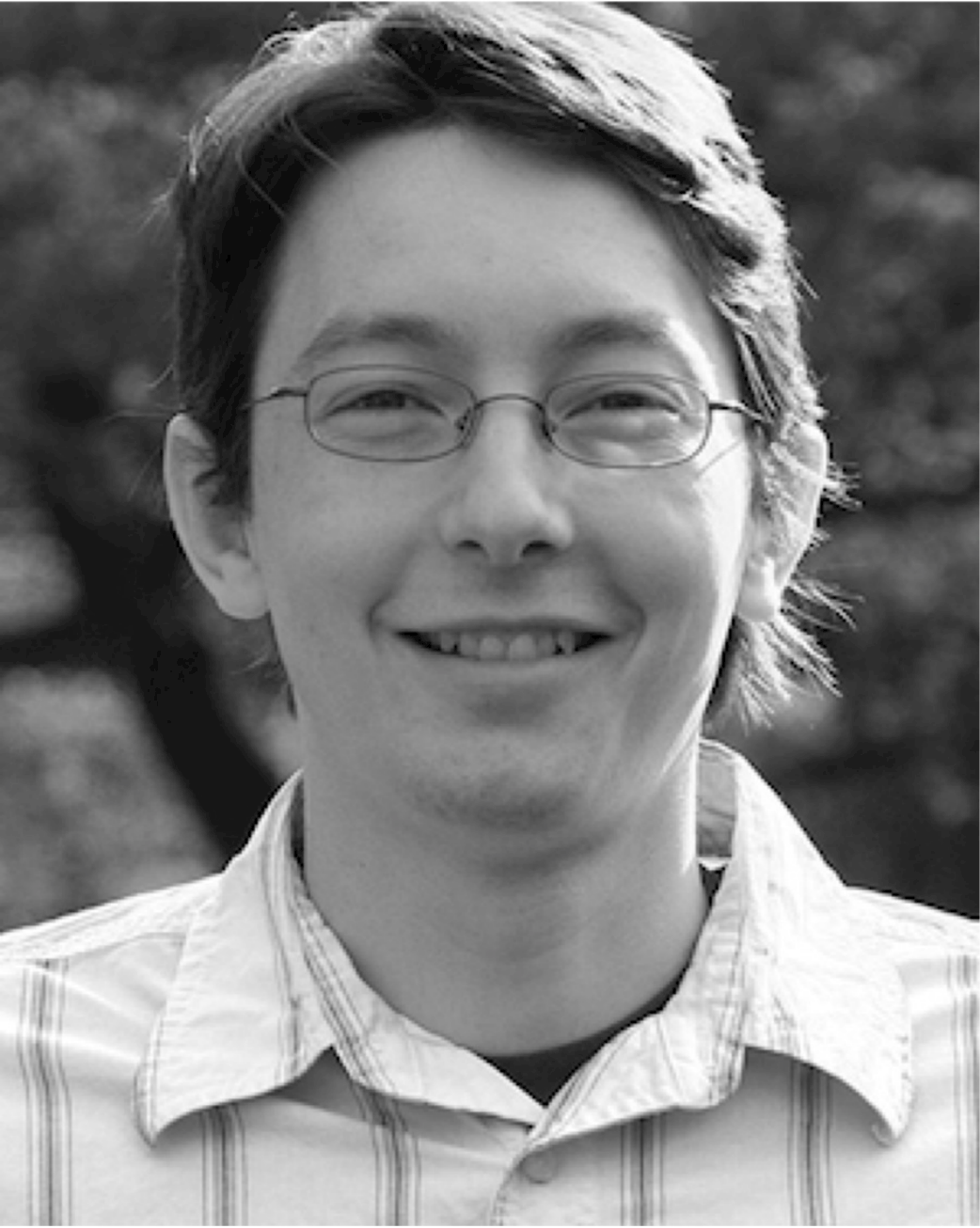}
  \end{center}  
  \vspace{-.215in}
\end{wrapfigure}
\noindent\textbf{Steven L. Brunton} received the B.S. degree in mathematics with a minor in control and dynamical systems from the California Institute of Technology, Pasadena, CA, in 2006, and the Ph.D. degree in mechanical and aerospace engineering from Princeton University, Princeton, NJ, in 2012.  He is currently an Acting Assistant Professor of applied mathematics at the University of Washington.\\

\begin{wrapfigure}{l}{0.15\textwidth}
\vspace{-.15in}
  \begin{center}
    \includegraphics[width=.15\textwidth]{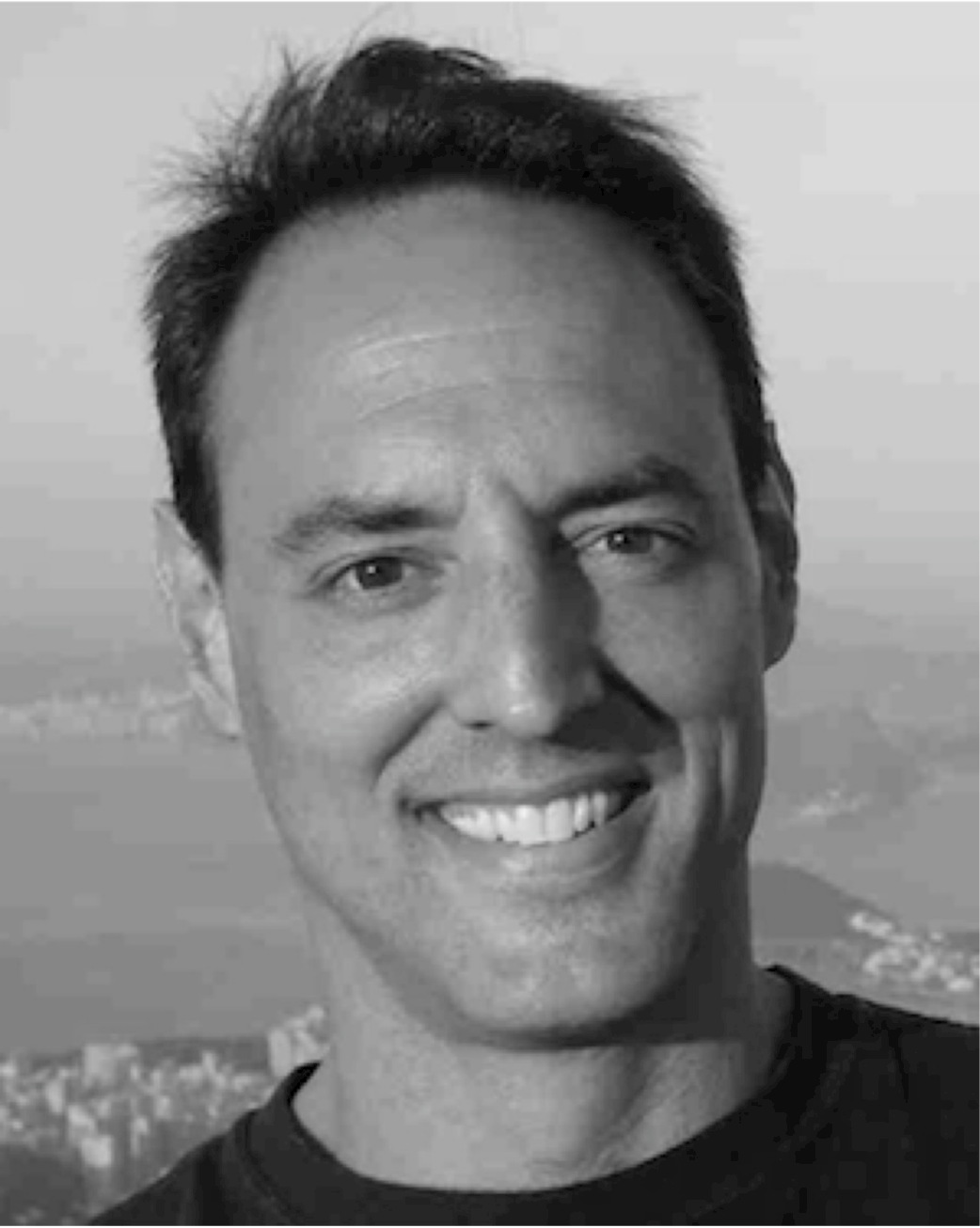}
  \end{center}
  \vspace{-.275in}
\end{wrapfigure}
\noindent\textbf{J. Nathan Kutz} received the B.S. degrees in physics and mathematics from the University of Washington, Seattle, WA, in 1990, and the Ph.D. degree in applied mathematics from Northwestern University, Evanston, IL, in 1994.  He is currently a Professor and Chair of applied mathematics at the University of Washington.

\end{document}